\documentclass[11pt,a4paper]{article}

\usepackage{theorem,enumerate}
\usepackage{amsmath,latexsym,amssymb,amsfonts}
\usepackage{mathrsfs}
\usepackage[dvipdf]{graphicx}

\newcounter{Scounter}
\theorembodyfont{\normalfont\slshape}
\newtheorem{thm}{Theorem}

\newtheorem{Thm}{Theorem}

\newtheorem{lem}{Lemma}

\newtheorem{Q}{Question}

\newtheorem{Claim}{Claim}[section]
\newtheorem{Fact}[Claim]{Fact}

\newtheorem{prob}{\normalfont\itshape Problem}
\newtheorem{remark}{\normalfont\itshape Remark}

\newcommand{\proof}{\medbreak\noindent\textit{Proof.}\quad}

\newcommand{\qed}{{$\quad\square$\vs{3.6}}}

\newcommand{\bss}{\setminus}

\newcommand{\vs}[1]{\vspace*{#1 mm}}

\newcommand{\ora}{\overrightarrow}
\newcommand{\ola}{\overleftarrow}

\numberwithin{equation}{section}

\def\Vec#1{\mbox{\boldmath $#1$}}

\def\calH{{ \mathcal{H}}}

\makeatletter
\def\thanks#1{%
   \footnotemark
   \edef\@tempa{\noexpand\noexpand\noexpand\footnotetext[\the\c@footnote]}%
   \toks@\expandafter{\@thanks}%
   \toks\tw@{{#1}}
   \xdef\@thanks{\the\toks@\@tempa\the\toks\tw@}}
\makeatother

\bfseries\normalfont

\begin{document}

\title{Dominating cycles and forbidden pairs \\containing a path of order $5$}

\author{
Shuya Chiba$^{1}$\thanks{This work was partially supported by JSPS KAKENHI grant 26800083\\
\hspace{+14pt}
E-mail address: \texttt{schiba@kumamoto-u.ac.jp}}
\and
Michitaka Furuya$^{2}$\thanks{E-mail address: \texttt{michitaka.furuya@gmail.com} } 
\and 
Shoichi Tsuchiya$^{3}$\footnote{E-mail address: \texttt{wco.liew6.te@gmail.com}}
\vspace{+8pt}
\\
$^1$\small\textsl{Department of Mathematics and Engineering, 
Kumamoto University}\\ 
\small\textsl{2-39-1, Kurokami, Kumamoto 860-8555, Japan}\\
\small
$^{2}$\small\textsl{Department of Mathematical Information Science,
Tokyo University of Science,}\\ 
\small\textsl{1-3 Kagurazaka, Shinjuku-ku, Tokyo 162-8601, Japan}\\
\small
$^{3}$\small\textsl{Department of Mathematics,
Keio University,}\\ 
\small\textsl{3-14-1 Hiyoshi, Kohoku-ku, Yokohama-shi, Kanagawa 223-8522, Japan}\\
}

\date{}

\maketitle

\vspace{-20pt}
\begin{abstract}
A cycle $C$ in a graph $G$ is dominating 
if every edge of $G$ is incident with a vertex of $C$. 
For a set $\calH$ of connected graphs, 
a graph $G$ is said to be $\calH$-free 
if $G$ does not contain any member of $\calH$ as an induced subgraph. 
When $|\calH| = 2$, $\calH$ is called a forbidden pair. 
In this paper, 
we investigate the characterization of the class of the forbidden pairs guaranteeing the existence of a dominating cycle 
and show the following two results: 
(i) Every $2$-connected $\{P_{5}, K_{4}^{-}\}$-free graph contains a longest cycle which is a dominating cycle. 
(ii) Every $2$-connected $\{P_{5}, W^{*}\}$-free graph contains a longest cycle which is a dominating cycle. 
Here $P_{5}$ is the path of order $5$, 
$K_{4}^{-}$ is the graph obtained from the complete graph of order $4$ by removing one edge, 
and $W^{*}$ is a graph obtained from two triangles and an edge 
by identifying one vertex in each.

\end{abstract}

\section{Introduction}
\label{introduction}

In this paper, we consider only finite simple graphs. 
For terminology and notation not defined in this paper, we refer the readers to \cite{D}. 
A graph $G$ is said to be \textit{Hamiltonian} if $G$ has a \textit{Hamilton cycle}, i.e., a cycle containing all vertices of $G$.
A cycle $C$ in a graph $G$ is \textit{dominating} 
if every edge of $G$ is incident with a vertex of $C$.

Let $\calH$ be a set of connected graphs. 
A graph $G$ is said to be \textit{$\calH$-free} 
if $G$ does not contain $H$ as an induced subgraph for all $H$ in $\calH$, 
and we call each graph $H$ of $\calH$ a \textit{forbidden subgraph}. 
We call $\calH$ a \textit{forbidden pair} if $|\calH| = 2$. 
When we consider $\calH$-free graphs, 
we assume that each member of $\calH$ has order at least $3$ 
because $K_{2}$ is the only connected graph of order $2$ 
and $K_{1}$ is the unique $K_{2}$-free connected graph
(here $K_{n}$ denotes the complete graph of order $n$). 
In order to state results clearly, 
we further introduce the following notation. 
For two sets $\calH_{1}$ and $\calH_{2}$ of connected graphs, 
we write $\calH_{1} \le \calH_{2}$ 
if for every graph $H_{2}$ in $\calH_{2}$, 
there exists a graph $H_{1}$ in $\calH_{1}$ 
such that $H_{1}$ is an induced subgraph of $H_{2}$. 
Note that if $\calH_{1} \le \calH_{2}$, then every $\calH_{1}$-free graph is also $\calH_{2}$-free.

\begin{figure}[htbp]
\begin{center}
\hspace{-10pt}\includegraphics[scale=0.9,clip]{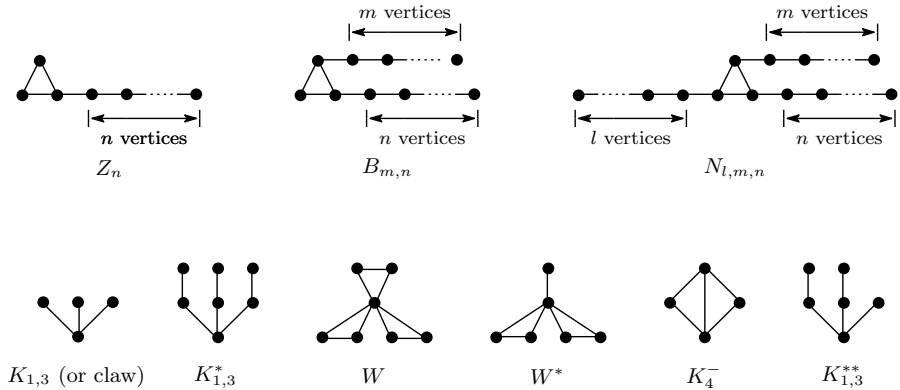}
\vspace{-25pt}
\caption{Forbidden subgraphs}
\label{forbidden(Fig)}
\end{center}
\end{figure}

The forbidden pairs that force the existence of a Hamilton cycle in $2$-connected graphs had been studied in \cite{BV,DGJ,GJ}.
In 1991, a characterization of such pairs was accomplished by Bedrossian \cite{Be}. 
Later, Faudree and Gould \cite{FG} extended the result of Bedrossian 
by regarding finite number of $2$-connected $\{H_{1}, H_{2}\}$-free non-Hamiltonian graphs as exceptions. 
Here let $P_{n}$ denote the path of order $n$, 
and the graphs $K_{1, 3}$ (or claw), $Z_{n}$, 
$B_{m, n}$ and $N_{l, m, n}$ 
are the ones that are depicted in Figure~\ref{forbidden(Fig)}.

\begin{Thm}[Faudree and Gould \cite{FG}]
\label{FG}
Let $\calH$ be a forbidden pair. 
Then every $2$-connected $\calH$-free graph of sufficiently large order is Hamiltonian 
if and only if 
$\calH \le \{K_{1, 3}, P_{6}\}$, 
$\calH \le \{K_{1, 3}, Z_{3}\}$, 
$\calH \le \{K_{1, 3}, B_{1, 2}\}$, 
or 
$\calH \le \{K_{1, 3}, N_{1, 1, 1}\}$. 
\end{Thm}

The purpose of this paper 
is to consider the analogue of Theorem \ref{FG} for dominating cycles which are relaxed structures of a Hamilton cycle. 
More precisely, 
we consider the following problem.

\begin{prob}
\label{determine H}
Determine the set $\Vec{H}$ (resp., $\Vec{H}'$) of forbidden pairs $\calH$ which satisfy that
every $2$-connected $\calH$-free graph 
(resp., every $2$-connected $\calH$-free graph of sufficiently large order) has a dominating cycle.
\end{prob}

Concerning this problem, 
the authors 
proved the following result in \cite{CFT} 
(here let $K_{1, 3}^{*}$, $W$, $W^{*}$ and $K_{4}^{-}$ 
be the ones that are depicted 
in Figure~\ref{forbidden(Fig)}).

\begin{Thm}[\cite{CFT}]
\label{2-conn (necessity)}
Let $\cal{H}$ be a forbidden pair. 
If 
there exists a positive integer $n_{0} = n_{0}(\calH)$ such that 
every $2$-connected $\cal{H}$-free graph of order at least $n_{0}$ has a dominating cycle, 
then 
$\calH \le \{K_{1, 3}, Z_{4}\}$, 
$\calH \le \{K_{1, 3}, B_{1, 2}\}$, 
$\calH \le \{K_{1, 3}, N_{1, 1, 1}\}$, 
$\calH \le \{ P_{4}, W \}$, 
$\calH \le \{ K_{1, 3}^{*}, Z_{1} \}$, 
$\calH \le \{ P_{5}, W^{*} \}$, 
or 
$\calH \le \{ P_{5}, K_{4}^{-}\}$. 
\end{Thm}

In the same paper, 
the authors also conjectured that 
the converse of Theorem~\ref{2-conn (necessity)} holds 
and gave a partial solution of the conjecture as follows. 
Here $K_{1, 3}^{**}$ is the graph obtained from $K_{1, 3}^{*}$ by deleting one leaf 
(see Figure~\ref{forbidden(Fig)}).

\begin{Thm}[\cite{CFT}]
\label{2-conn (sufficiency)}
If $\calH \le \{K_{1, 3}, Z_{4}\}$, 
$\calH \le \{K_{1, 3}, B_{1, 2}\}$, 
$\calH \le \{K_{1, 3}, N_{1, 1, 1}\}$, 
$\calH \le \{ P_{4}, W \}$, 
or 
$\calH \le \{  K_{1, 3}^{**}, Z_{1} \}$, 
then 
every $2$-connected $\cal{H}$-free graph has a dominating cycle. 
\end{Thm}

In this paper, 
we show that the above conjecture is also true for the cases where 
$\calH \le \{ P_{5}, W^{*} \}$ 
and 
$\calH \le \{ P_{5}, K_{4}^{-}\}$ 
by considering slightly stronger statements.

\begin{thm}
\label{P_{5}, W^{*}}
Every $2$-connected $\{ P_{5}, W^{*}  \}$-free graph contains a longest cycle which is a dominating cycle.
\end{thm}

\begin{thm}
\label{P_{5}, K_{4}^{-}}
Every $2$-connected $\{ P_{5}, K_{4}^{-} \}$-free graph contains a longest cycle which is a dominating cycle.
\end{thm}

\begin{remark}
By Theorems \ref{2-conn (necessity)}, \ref{2-conn (sufficiency)}, \ref{P_{5}, W^{*}} and \ref{P_{5}, K_{4}^{-}}, 
the remaining problem is only that whether the pair $\{K_{1, 3}^{*}, Z_{1}\}$ belongs to the class $\Vec{H}$ (resp., $\Vec{H}'$) 
of Problem \ref{determine H} 
or not.
Olariu \cite{O} showed that if a connected $Z_{1}$-free graph $G$ contains a triangle, then $G$ is a complete multipartite graph.
On the other hand, it is easy to check that every $2$-connected complete multipartite graph containing a triangle has a dominating cycle.
Thus the pair $\{K_{1,3}^{*},Z_{1}\}$ belongs to the class $\Vec{H}$ (resp., $\Vec{H}'$) if and only if 
the pair $\{K_{1,3}^{*},K_{3}\}$ belongs to the class $\Vec{H}$ (resp., $\Vec{H}'$).
Consequently, we can deduce the target pair to $\{K_{1,3}^{*},K_{3}\}$.
Although we do not know the answer at the moment, 
we believe that the pair $\{K_{1, 3}^{*}, K_{3}\}$ belongs to the class.
\end{remark}

\medskip

In Section \ref{preparation}, 
we will introduce the lemmas in order to show Theorems \ref{P_{5}, W^{*}} and \ref{P_{5}, K_{4}^{-}}, 
and we prove Theorems \ref{P_{5}, W^{*}} and \ref{P_{5}, K_{4}^{-}} in Sections \ref{proof of P_{5},W^{*}} and \ref{proof of P_{5},K_{4}^{-}}, 
respectively.

\section{Preparation for the proofs of Theorems \ref{P_{5}, W^{*}} and \ref{P_{5}, K_{4}^{-}}}
\label{preparation}

In this section, we prepare lemmas which will be used in the proofs of Theorems \ref{P_{5}, W^{*}} and \ref{P_{5}, K_{4}^{-}}.
To do that, 
we first prepare 
terminology and notation which we use in the rest. 

Let $G$ be a graph. 
We denote by $V(G)$ and $E(G)$ 
the vertex set and the edge set of $G$, respectively, 
and let $|G| = |V(G)|$. 
For $X \subseteq V(G)$, 
we let $G[X]$ denote the subgraph induced by $X$ in $G$, 
and let $G - X = G[V(G) \setminus X]$. 
Let $v$ be a vertex of $G$. 
We denote by $N_{G}(v)$ 
the neighborhood 
of $v$ in $G$. 
For $X \subseteq V(G) \bss \{v\}$, 
we let $N_{G}(v; X) = N_{G}(v) \cap X$, 
and for $V, X \subseteq V(G)$ with $V \cap X = \emptyset$, 
let $N_{G}(V; X) = \bigcup_{v \in V}N_{G}(v; X)$. 
In this paper, 
we often identify a subgraph $F$ of $G$ with its vertex set $V(F)$
(for example, $N_{G}(v; V(F))$ is often denoted by $N_{G}(v; F)$).

A path with ends $u$ and $v$ 
is denoted by a \textit{$(u, v)$-path}. 
For a subgraph $H$ of $G$, 
a path $P$ of $G$ such that $|P| \ge 2$ is called a $H$-\textit{path} 
if ends of $P$ only belong to $H$. 
We write a cycle (or a path) $C$ with a given orientation by $\ora{C}$.
If there exists no chance of confusion, we abbreviate $\ora{C}$ by $C$. 
Let $\ora{C}$ be an oriented cycle or a path.
For $u, v \in V(C)$,
we denote by $u\ora{C}v$
the $(u, v)$-path on $\ora{C}$.
The reverse sequence
of $u\ora{C}v$
is denoted by $v\ola{C}u$.
For $v \in V(C)$, 
we denote the $h$-th successor and the $h$-th predecessor of $v$ 
on $\ora{C}$ by $v^{+h}$ and $v^{-h}$, respectively, 
and let $v^{+0} = v^{-0} = v$. 
For $X \subseteq V(C)$, 
we define $X^{+h} = \{v^{+h} : v \in X \}$ 
and $X^{-h} = \{v^{-h} : v \in X\}$, respectively.  
We abbreviate $v^{+1}$, $v^{-1}$, $X^{+1}$ and $X^{-1}$ by $v^{+}$, $v^{-}$, $X^+$ and $X^-$, respectively.

\subsection{Lemmas for $P_{5}$-free graphs}
\label{basic forbid}

In this subsection, 
we give the following two lemmas (Lemmas \ref{consecutive} and \ref{if V(H) - N_{G}(v_{i}) neq emptyset}) 
to make it easy to use the assumption ``$P_{5}$-free'' 
in the proofs of Theorems \ref{P_{5}, W^{*}} and \ref{P_{5}, K_{4}^{-}}.

\begin{lem}
\label{consecutive} 
Let $G$ be a graph, 
and let $Q_{1}$ and $Q_{2}$ be paths of order at least $3$ with a common end $a$ 
such that $Q_{1} - a$ and $Q_{2} - a$ are vertex-disjoint. 
If $G$ is $P_{5}$-free and $Q_{1}$ is an induced path, 
then 
$N_{G}(Q_{1} - a; Q_{2} - a) \neq \emptyset$ 
or 
$V(Q_{2}) \setminus \{a\} \subseteq N_{G}(a)$. 
\end{lem}
\noindent
\textit{Proof of Lemma \ref{consecutive}.}~
Suppose that $N_{G}(Q_{1} - a; Q_{2} - a) = \emptyset$ and $V(Q_{2}) \setminus \{a\} \not\subseteq N_{G}(a)$.
Write $Q_{1}=a_{1}a_{2} \dots a_{l}$ and $Q_{2}=a'_{1}a'_{2} \dots a'_{l'}$, where $a_{1}=a'_{1}=a$.
Let $i~(1\leq i\leq l')$ be the minimum index with $aa'_{i}\not\in E(G)$.
Note that $a\not=a'_{i-1}$ and $aa'_{i-1}\in E(G)$.
Hence $a_{3}a_{2}aa'_{i-1}a'_{i}$ is an induced path of $G$ 
because $Q_{1}$ is an induced path and $N_{G}(Q_{1} - a; Q_{2} - a) = \emptyset$, 
which is a contradiction.
\qed

By Lemma \ref{consecutive}, 
we can easily obtain the following.

\begin{lem}
\label{if V(H) - N_{G}(v_{i}) neq emptyset}
Let $G$ be a $P_{5}$-free graph, $\ora{C}$ be a cycle and $H$ be a component of $G - C$, 
and let $v \in N_{G}(H; C)$ such that $V(H) \setminus N_{G}(v) \neq \emptyset$.
If $N_{G}(H; v^{+} \ora{C} a) = \emptyset$ for some $a \in V(C) \setminus \{v, v^{+}\}$ 
(resp. $N_{G}(H; v^{-} \ola{C} a) = \emptyset$ for some $a \in V(C) \setminus \{v, v^{-}\}$), 
then 
$V(v^{+} \ora{C} a) \subseteq N_{G}(v)$ 
(resp. $V(v^{-} \ola{C} a) \subseteq N_{G}(v)$). 
\end{lem}
\noindent
\textit{Proof of Lemma \ref{if V(H) - N_{G}(v_{i}) neq emptyset}.}~
By the symmetry, it suffice to consider the case where $N_{G}(H; v^{+} \ora{C} a) = \emptyset$ for some $a \in V(C) \setminus \{v, v^{+}\}$.
Since $V(H) \setminus N_{G}(v) \neq \emptyset$, there exist two vertices $u,u'\in V(H)$ such that $vu,uu'\in E(G)$ and $vu'\not\in E(G)$.
Now we take two paths $Q_{1}=vuu'$ and $Q_{2}=v\ora{C}a$.
Then $Q_{1}$ is an induced path of $G$ and $N_{G}(Q_{1} - v; Q_{2} - v) = \emptyset$.
This together with Lemma \ref{consecutive} leads to $V(v^{+} \ora{C} a)=V(Q_{2}) \setminus \{v\} \subseteq N_{G}(v)$.
\qed

\subsection{Properties of longest cycles in graphs}
\label{basic}

In this subsection, 
we introduce the basic lemmas concerning the properties of longest cycles in graphs.

\begin{figure}[htbp]
\begin{center}
\hspace{-10pt}\includegraphics[scale=0.9,clip]{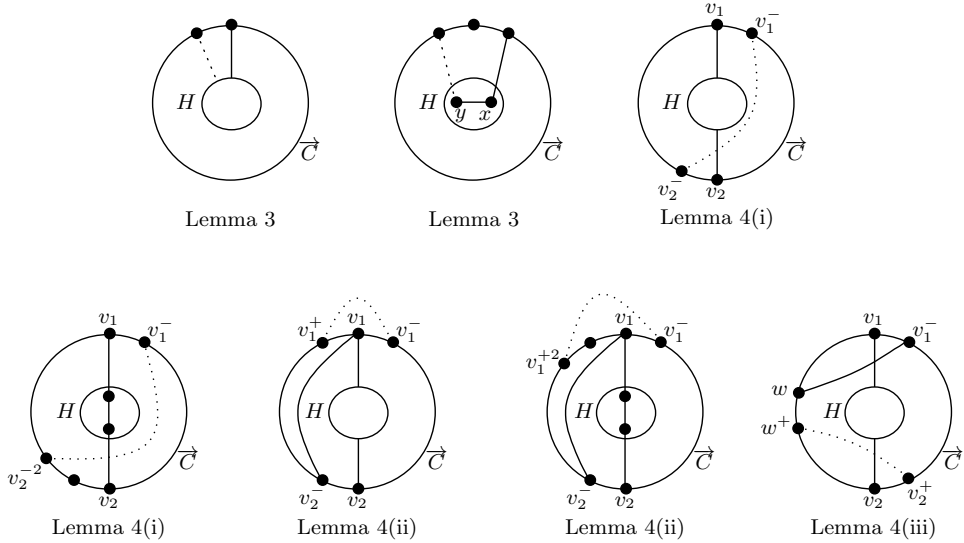}
\vspace{-25pt}
\caption{Longest cycles in graphs}
\label{longestlemma(Fig)}
\end{center}
\end{figure}

We fix the following notation in this subsection. 
Let $G$ be a graph and $\ora{C}$ be a longest cycle of $G$, and let $H$ be a component of $G - C$. 
Then the following two lemmas hold (Lemmas \ref{maximality of |V(C)| (1)} and \ref{maximality of |V(C)| (2)}). 
Since the proofs directly follow from the maximality of $|C|$, 
we omit it (see also Figure \ref{longestlemma(Fig)}).

\begin{lem}
\label{maximality of |V(C)| (1)}
$N_{G}(x; C) \cap N_{G}(y; C)^{-} = \emptyset$ for $x, y \in V(H)$. 
In particular, 
if $x \neq y$, then 
$N_{G}(x; C) \cap N_{G}(y; C)^{-2} = \emptyset$. 
\end{lem}

\begin{lem}
\label{maximality of |V(C)| (2)}
Let $v_{1}$ and $v_{2}$ be two distinct vertices in $N_{G}(H; C)$. 
Then the following hold. 
\begin{enumerate}[{\upshape(i)}]

\item
\label{no v1-v2-} 
There exists no $C$-path joining $v_{1}^{-}$ and $v_{2}^{-}$, 
and joining $v_{1}^{+}$ and $v_{2}^{+}$, respectively; 
in particular, 
$E(G) \cap \{v_{1}^{-}v_{2}^{-}, v_{1}^{+}v_{2}^{+}\} = \emptyset$. 
Moreover, 
if $|N_{G}(v_{1}; H) \cup N_{G}(v_{2}; H)| \ge 2$, then 
there exists no $C$-path joining $v_{1}^{-}$ and $v_{2}^{-2}$, 
and joining $v_{1}^{+}$ and $v_{2}^{+2}$, respectively; 
in particular, 
$E(G) \cap \{v_{1}^{-}v_{2}^{-2}, v_{1}^{+}v_{2}^{+2}\} = \emptyset$.

\item
\label{no v1-w v2-w-2} 
If $v_{1}v_{2}^{-} \in E(G)$, then $v_{1}^{-}v_{1}^{+} \notin E(G)$. 
Moreover, 
if $v_{1}v_{2}^{-} \in E(G)$
and 
$|N_{G}(v_{1}; H) \cup N_{G}(v_{2}; H)| \ge 2$, 
then $E(G) \cap \{v_{1}^{-}v_{1}^{+2}, v_{1}^{-2}v_{1}^{+}\} = \emptyset$.

\item
\label{no v1-w v2+w+} 
If $v_{1}^{-}w \in E(G)$ for some vertex $w$ in $v_{1} \ora{C}v_{2}^{-}$, then $v_{2}^{+}w^{+} \notin E(G)$. 
If $v_{1}^{+}w \in E(G)$ for some vertex $w$ in $v_{1} \ola{C}v_{2}^{+}$, then $v_{2}^{-}w^{-} \notin E(G)$. 

\end{enumerate}
\end{lem}

\subsection{Longest cycles in $P_{5}$-free graphs having no dominating longest cycle}
\label{longest cycle P5-free no dom}

For a cycle $C$ of a graph $G$, let 
$\mu(C) 
= \max \{ |F| : F$ is a component of $G - C \}$, 
and we define $\omega(C) = | \{ F : F$ is a component of $G - C$ such that $|F| = \mu(C) \}|$.

Now let $G$ be a graph, 
and we suppose that any longest cycles of $G$ are not dominating cycles 
(i.e., $\mu(C) \ge 2$ for every longest cycle $C$ of $G$),  
and 
let $\ora{C}$ be a longest cycle of $G$. 
Suppose further that $C$ was chosen so that 
\begin{enumerate}
\newcounter{mymemory}
\renewcommand{\labelenumi}{\upshape{(C\arabic{enumi}})}
\item
\label{choice 1}
$\mu(C)$ is as small as possible, and  
\vspace{-6pt} 

\item 
\label{choice 2}
$\omega(C)$ is as small as possible, subject to (C\ref{choice 1}). 

\setcounter{mymemory}{\value{enumi}}
\end{enumerate}
Let $H$ be a component of $G-C$ such that 
$|H| = \mu(C) \ ( \ \ge 2)$.

\begin{lem}
\label{indep. set S}
If $S$ is an independent set of $G$ such that $S \subseteq V(C)$ and $N_{G}(S; G- C) = \emptyset$, 
then there exists no longest cycle $D$ of $G$ such that $V(D) \supseteq V(C) \bss S$ and $V(D) \cap V(H) \neq \emptyset$. 
\end{lem}
\noindent
\textit{Proof of Lemma \ref{indep. set S}.}~Suppose that there exists a longest cycle $D$ of $G$ 
such that $V(D) \supseteq V(C) \bss S$ and $V(D) \cap V(H) \neq \emptyset$. 
Let $H'$ be an arbitrary component of $G - D$. 
By the assumptions of $S$, we see that $|H'| = 1$ or $H'$ is an induced subgraph of some component of $G - C$
because $V(C) \setminus S \subseteq V(D)$. 
This implies that $\mu(D) \le \mu(C)$. 
Moreover, if $\mu(D) = \mu(C)$, then 
the number of components of order $\mu(C)$ in $G - D$ is less than $\omega(C)$ 
because $V(D) \cap V(H) \neq \emptyset$. 
This contradicts the choice (C1) or (C2). 
\qed

By Lemma \ref{indep. set S}, 
the following two lemmas hold for $P_{5}$-free graphs.

\begin{figure}[htbp]
\begin{center}
\hspace{-10pt}\includegraphics[scale=0.9,clip]{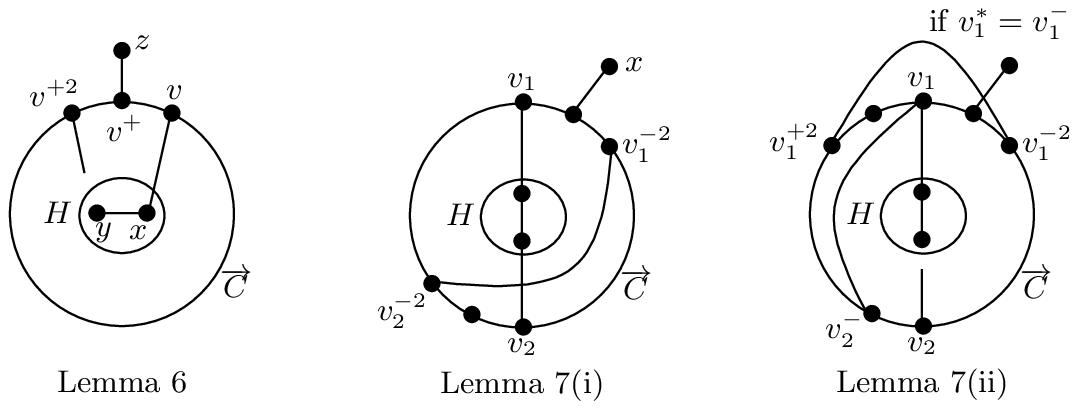}
\vspace{-25pt}
\caption{Lemmas \ref{N_{G}(H; C) cap N_{G}(H; C)^{-2} = emptyset} and \ref{v1 v2 in P5-free}}
\label{nodomlongestlemma(Fig)}
\end{center}
\end{figure}

\begin{lem}
\label{N_{G}(H; C) cap N_{G}(H; C)^{-2} = emptyset}
If $G$ is $P_{5}$-free, 
then $N_{G}(H; C) \cap N_{G}(H; C)^{-2} = \emptyset$. 
\end{lem}
\noindent
\textit{Proof of Lemma \ref{N_{G}(H; C) cap N_{G}(H; C)^{-2} = emptyset}.}~Suppose that 
$N_{G}(H; C) \cap N_{G}(H; C)^{-2} \neq \emptyset$, 
and let $v \in N_{G}(H; C) \cap N_{G}(H; C)^{-2}$. 
Note that by Lemma \ref{maximality of |V(C)| (1)}, 
$N_{G}(v^{+}; H) = \emptyset$. 
Note also that 
$G$ contains a longest cycle $D$ 
such that $V(D) \supseteq V(C) \bss \{v^{+}\}$ and $V(D) \cap V(H) \neq \emptyset$ 
because $v, v^{+2} \in N_{G}(H; C)$. 
Hence by Lemma \ref{indep. set S}, 
there exists a component $H'$ of $G - C$ such that 
$H' \neq H$ and $N_{G}(v^{+}; H') \neq \emptyset$. 
Let $x\in N_{G}(v;H)$, $y\in N_{H}(x)$ and 
$z \in N_{G}(v^{+}; H')$ (see Figure \ref{nodomlongestlemma(Fig)}). 
Consider the paths $Q_{1} = vxy$ and $Q_{2} = vv^{+}z$. 
Since $v, v^{+2} \in N_{G}(H; C)$, 
it follows from Lemma \ref{maximality of |V(C)| (1)} that $vy \notin E(G)$, 
and thus $Q_{1}$ is an induced path. 
Hence by Lemma \ref{consecutive}, 
$N_{G}(\{x, y\}; \{v^{+}, z\}) = N_{G}(Q_{1} - v; Q_{2} - v) \neq \emptyset$ 
or $zv \in E(G)$. 
Since $N_{G}(v^{+}; H) = N_{G}(H; H') = \emptyset$, 
we have 
$zv \in E(G)$, 
but this contradicts Lemma \ref{maximality of |V(C)| (1)}. 
\qed

\begin{lem}
\label{v1 v2 in P5-free}
Let $v_{1}$ and $v_{2}$ be two distinct vertices 
in $N_{G}(H; C)$ such that $|N_{G}(v_{1}; H) \cup N_{G}(v_{2}; H)| \ge 2$, 
and suppose that $G$ is $P_{5}$-free. 
Then the following hold. 
\begin{enumerate}[{\upshape(i)}]
\item 
\label{no v1-2v2-2}
$v_{1}^{-2}v_{2}^{-2} \notin E(G)$. 

\item 
\label{no v_{1}^{-2}v_{1}^{+2} if v_{1}v_{2}^{-} in E(G)}
If $V(H) \setminus N_{G}(v_{1}) \neq \emptyset$ and 
$v_{1}v_{2}^{-} \in E(G)$, then $v_{1}^{-2}v_{1}^{+2} \notin E(G)$. 
\end{enumerate}
\end{lem}
 
\noindent
\textit{Proof of Lemma \ref{v1 v2 in P5-free}.}~
Note that $v_{1}$, $v_{1}^{-}$, $v_{1}^{-2}$, $v_{2}$, $v_{2}^{-}$ and $v_{2}^{-2}$ are distinct vertices by Lemma \ref{maximality of |V(C)| (1)}.
Let $\ora{P}$ be a $(v_{1}, v_{2})$-path 
such that $|P| \ge 4$ and $V(P) \bss \{v_{1}, v_{2}\} \subseteq V(H)$. 

To show (\ref{no v1-2v2-2}), 
suppose that 
$v_{1}^{-2}v_{2}^{-2} \in E(G)$, 
and let $\ora{D} = v_{1}^{-2}v_{2}^{-2} \ola{C} v_{1} \ora{P} v_{2} \ora{C} v_{1}^{-2}$. 
Then $D$ is a cycle in $G$ such that $V(D) = (V(C) \bss \{v_{1}^{-}, v_{2}^{-}\}) \cup V(P)$. 
Hence by the maximality of $|C|$, $|P| = 4$. 
Since $\{v_{1}^{-}, v_{2}^{-}\}$ is an independent set of $G$ by Lemma \ref{maximality of |V(C)| (2)}(\ref{no v1-v2-}) 
and since $D$ is  also a longest cycle of $G$, 
it follows from Lemma \ref{indep. set S} that $N_{G}(v_{i}^{-}; G-C) \neq \emptyset$ for some $i$ with $i \in \{1, 2\}$. 
Suppose that $N_{G}(v_{1}^{-}; G-C) \neq \emptyset$, 
and let $x \in N_{G}(v_{1}^{-}; G-C)$. 
Note that by Lemma \ref{maximality of |V(C)| (1)}, $x \notin V(H)$ (see Figure \ref{nodomlongestlemma(Fig)}). 
Consider the paths $Q_{1} = v_{1}^{-2}v_{1}^{-}x$ and $Q_{2} = v_{1}^{-2}v_{2}^{-2}v_{2}^{-}$. 
By Lemma \ref{maximality of |V(C)| (1)}, 
$Q_{1}$ is an induced path. 
By Lemma \ref{maximality of |V(C)| (2)}(\ref{no v1-v2-}), 
$Q_{2}$ is also an induced path. 
Hence by Lemma \ref{consecutive}, 
$N_{G}(\{v_{1}^{-}, x\}; \{v_{2}^{-}, v_{2}^{-2}\}) \neq \emptyset$, 
but this contradicts Lemma \ref{maximality of |V(C)| (2)}(\ref{no v1-v2-}). 
Thus $N_{G}(v_{1}^{-}; G-C) = \emptyset$. 
By the symmetry of $v_{1}$ and $v_{2}$, 
we can get a contradiction for the case where $N_{G}(v_{2}^{-}; G-C) \neq \emptyset$. 
Thus (\ref{no v1-2v2-2}) holds.

To show (\ref{no v_{1}^{-2}v_{1}^{+2} if v_{1}v_{2}^{-} in E(G)}), suppose next that 
$V(H) \setminus N_{G}(v_{1}) \neq \emptyset$ and $\{v_{1}v_{2}^{-}, v_{1}^{-2}v_{1}^{+2}\} \subseteq E(G)$, 
and let $\ora{D'} = v_{2}^{-}v_{1} \ora{P} v_{2} \ora{C}$ $v_{1}^{-2} v_{1}^{+2} \ora{C} v_{2}^{-}$. 
Then $D'$ is a cycle in $G$ such that $V(D) = (V(C) \bss \{v_{1}^{-}, v_{2}^{-}\}) \cup V(P)$, 
and the maximality of $|C|$ implies that $|P| = 4$. 
Since $\{v_{1}^{-}, v_{1}^{+}\}$ is an independent set of $G$ by Lemma \ref{maximality of |V(C)| (2)}(\ref{no v1-w v2-w-2}) 
and since $D'$ is a longest cycle of $G$, 
it follows from Lemmas \ref{maximality of |V(C)| (1)} and  \ref{indep. set S} that 
$N_{G}(v_{1}^{*}; H') \neq \emptyset$ 
for some $v_{1}^{*} \in \{v_{1}^{-}, v_{1}^{+}\}$ 
and some component $H'$ of $G - C$ with $H' \neq H$ (see Figure \ref{nodomlongestlemma(Fig)}). 
Since $V(H) \setminus N_{G}(v_{1}) \neq \emptyset$, 
$G[V(H) \cup \{v_{1}\}]$ contains an induced path $Q_{1}'$ of order at least $3$ with an end $v_{1}$. 
By Lemma \ref{maximality of |V(C)| (1)}, 
$G[V(H') \cup \{v_{1}, v_{1}^{*}\}]$ 
contains an induced path $Q_{2}'$ of order at least $3$ with an end $v_{1}$ and $v_{1}v_{1}^{*} \in E(Q_{2}')$. 
Hence by Lemma \ref{consecutive} 
and since $N_{G}(H; H') = \emptyset$, 
we see that $N_{G}(v_{1}^{*}; H) \neq \emptyset$, 
which contradicts Lemma \ref{maximality of |V(C)| (1)}. 
Thus (\ref{no v_{1}^{-2}v_{1}^{+2} if v_{1}v_{2}^{-} in E(G)}) also holds. 
\qed

\section{Proof of Theorem \ref{P_{5}, W^{*}}}
\label{proof of P_{5},W^{*}}

Let $G$ be a $2$-connected $\{P_{5}, W^{*}\}$-free graph, 
and we show that $G$ contains a longest cycle which is a dominating cycle. 
By way of a contradiction, 
suppose that any longest cycles of $G$ are not dominating cycles.
Let $\ora{C}$ be the same described as in 
the paragraph preceding Lemma \ref{indep. set S} in Subsection \ref{longest cycle P5-free no dom}, 
and let $H$ be a component of $G-C$ such that $|H| = \mu(C) \ ( \ \ge 2)$. 
Since $G$ is $2$-connected, 
there exist two distinct vertices $v_{1}$ and $v_{2}$ in $N_{G}(H; C)$ 
such that $|N_{G}(v_{1}; H) \cup N_{G}(v_{2}; H)| \ge 2$. 
Then, $|v_{i} \ora{C} v_{3-i}| \ge 4$ for $i \in \{1, 2\}$ because $C$ is longest.
(Note that by these assumptions, in this proof, 
we can use all lemmas of Section \ref{preparation}.)  
We choose the vertices $v_{1}$ and $v_{2}$ so that 
\begin{align}
\label{N_{G}(H; v_{1}^{+} C v_{2}^{-}) is emptyset}
N_{G}(H; v_{1}^{+} \ora{C} v_{2}^{-}) = \emptyset. 
\end{align}

\begin{figure}[htbp]
\begin{center}
\hspace{-10pt}\includegraphics[scale=0.9,clip]{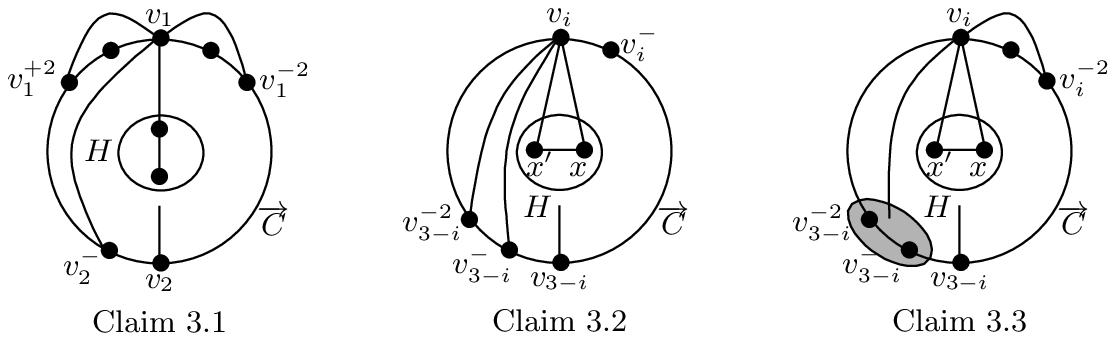}
\vspace{-25pt}
\caption{Claims \ref{xv1, xv2 in E(G)}--\ref{no v1v1-2if v1v2- or v1v2-2}}
\label{W(Fig)}
\end{center}
\end{figure}

\begin{Claim}
\label{xv1, xv2 in E(G)}
$V(H) \subseteq N_{G}(v_{1}) \cap N_{G}(v_{2})$. 
\end{Claim}
\proof 
Suppose that $V(H) \bss N_{G}(v_{1}) \neq \emptyset$. 
Then by Lemma \ref{if V(H) - N_{G}(v_{i}) neq emptyset} and (\ref{N_{G}(H; v_{1}^{+} C v_{2}^{-}) is emptyset}), 
we have $\{v_{1}v_{1}^{+2},$ $v_{1}v_{2}^{-}\} \subseteq E(G)$. 
Moreover, by Lemmas \ref{maximality of |V(C)| (1)} and \ref{N_{G}(H; C) cap N_{G}(H; C)^{-2} = emptyset},
$N_{G}(H;\{v_{1}^{-},v_{1}^{-2}\})=\emptyset $, and hence Lemma \ref{if V(H) - N_{G}(v_{i}) neq emptyset} yields that $v_{1}v_{1}^{-2}\in E(G)$.
Since $v_{1}v_{2}^{-} \in E(G)$, 
it follows from 
Lemma \ref{maximality of |V(C)| (2)}(\ref{no v1-w v2-w-2}) that 
$E(G) \cap \{v_{1}^{-}v_{1}^{+}, v_{1}^{-2}v_{1}^{+},$ $v_{1}^{-}v_{1}^{+2}\} = \emptyset$. 
By Lemma \ref{v1 v2 in P5-free}(\ref{no v_{1}^{-2}v_{1}^{+2} if v_{1}v_{2}^{-} in E(G)}), 
we also have $v_{1}^{-2}v_{1}^{+2} \notin E(G)$. 
Therefore, 
since $N_{G}(H; \{v_{1}^{-}, v_{1}^{-2}, v_{1}^{+},$ $v_{1}^{+2}\}) = \emptyset$ 
by Lemmas \ref{maximality of |V(C)| (1)} and \ref{N_{G}(H; C) cap N_{G}(H; C)^{-2} = emptyset}, 
we see that $G[\{v_{1}, v_{1}^{-}, v_{1}^{-2}, v_{1}^{+}, v_{1}^{+2}\} \cup N_{G}(v_{1}; H)]$ contains a $W^{*}$ as an induced subgraph 
(see Figure \ref{W(Fig)}), a contradiction. 
Thus $V(H) \subseteq N_{G}(v_{1})$. 
Similarly, we have $V(H) \subseteq N_{G}(v_{2})$. 
\qed

\begin{Claim}
\label{no v1v2- or v1v2-2}
For $i \in \{1, 2\}$, $|E(G) \cap \{v_{i}v_{3-i}^{-}, v_{i}v_{3-i}^{-2}\}| \le 1$. 
\end{Claim}
\proof 
Suppose that $\{v_{i}v_{3-i}^{-}, v_{i}v_{3-i}^{-2} \} \subseteq E(G)$, and let $xx' \in E(H)$. 
By Claim \ref{xv1, xv2 in E(G)}, $\{xv_{i}, x'v_{i}\} \subseteq E(G)$ (see Figure \ref{W(Fig)}).
By Lemmas \ref{maximality of |V(C)| (1)} and \ref{N_{G}(H; C) cap N_{G}(H; C)^{-2} = emptyset}, 
we have $E(G) \cap \{xv_{i}^{-},xv_{3-i}^{-}, xv_{3-i}^{-2}, x'v_{i}^{-}, x'v_{3-i}^{-}, x'v_{3-i}^{-2}\}$ $= \emptyset$. 
By Lemma \ref{maximality of |V(C)| (2)}(\ref{no v1-v2-}), 
we also have $E(G) \cap \{v_{i}^{-}v_{3-i}^{-}, v_{i}^{-}v_{3-i}^{-2}\} = \emptyset$. 
This implies that $G[\{v_{i}, v_{i}^{-}, v_{3-i}^{-},$ $v_{3-i}^{-2}, x, x'\}] \cong W^{*}$, a contradiction. 
Thus $|E(G) \cap \{v_{i}v_{3-i}^{-}, v_{i}v_{3-i}^{-2}\}| \le 1$. 
\qed

\begin{Claim}
\label{no v1v1-2if v1v2- or v1v2-2}
For $i \in \{1, 2\}$, 
if $E(G) \cap \{v_{i}v_{3-i}^{-}, v_{i}v_{3-i}^{-2}\} \neq \emptyset$, 
then $v_{i}v_{i}^{-2} \notin E(G)$. 
\end{Claim}
\proof 
Let $v^{*} \in \{v_{3-i}^{-}, v_{3-i}^{-2} \}$, 
and 
we show that if $v_{i}v^{*} \in E(G)$, then $v_{i}v_{i}^{-2} \notin E(G)$. 
By way of a contradiction, 
suppose that $\{v_{i}v^{*}, v_{i}v_{i}^{-2}\} \subseteq E(G)$, 
and let $xx' \in E(H)$. 
By Claim \ref{xv1, xv2 in E(G)}, $\{xv_{i}, x'v_{i}\} \subseteq E(G)$ (see Figure \ref{W(Fig)}).
By Lemmas \ref{maximality of |V(C)| (2)}(\ref{no v1-v2-}) and \ref{v1 v2 in P5-free}(\ref{no v1-2v2-2}), 
$E(G) \cap \{ v_{i}^{-}v^{*}, v_{i}^{-2}v^{*}\} = \emptyset$.  
By Lemmas \ref{maximality of |V(C)| (1)} and \ref{N_{G}(H; C) cap N_{G}(H; C)^{-2} = emptyset}, 
$E(G) \cap \{ xv_{i}^{-}, xv_{i}^{-2}, xv^{*}, x'v_{i}^{-}, x'v_{i}^{-2}, x'v^{*} \} = \emptyset$. 
This implies that 
$G[\{v_{i}, v_{i}^{-}, v_{i}^{-2}, x, x', v^{*} \}] \cong W^{*}$, a contradiction. 
Thus if $v_{i}v^{*} \in E(G)$, then $v_{i}v_{i}^{-2} \notin E(G)$. 
\qed

\begin{Claim}
\label{v1v2}
$v_{1}v_{2} \in E(G)$. 
\end{Claim}
\proof 
Let $x \in V(H)$, 
and 
consider the paths $Q_{1} = xv_{1}v_{1}^{-}$ and $Q_{2} = xv_{2}v_{2}^{-}$. 
By Lemma \ref{maximality of |V(C)| (1)}, each $Q_{i}$ is an induced path. 
Hence 
by Lemma \ref{consecutive}, 
$N_{G}(\{ v_{1}, v_{1}^{-} \}; \{ v_{2}, v_{2}^{-} \}) \neq \emptyset$. 
Since $v_{1}^{-}v_{2}^{-} \notin E(G)$ by Lemma \ref{maximality of |V(C)| (2)}(\ref{no v1-v2-}), 
we have 
$E(G) \cap \{v_{1}v_{2}^{-}, v_{1}^{-}v_{2}, v_{1}v_{2}\} \neq \emptyset$. 

Suppose that $v_{i}v_{3-i}^{-} \in E(G)$ for some $i\in \{1,2\}$, 
and next consider the paths $Q_{1}' = v_{i}v_{i}^{-}v_{i}^{-2}$ and $Q_{2}' = v_{i}v_{3-i}^{-}v_{3-i}^{-2}$. 
It follows from Claims \ref{no v1v2- or v1v2-2} and \ref{no v1v1-2if v1v2- or v1v2-2} that 
each $Q_{i}'$ is an induced path. 
Hence by Lemma \ref{consecutive}, 
$N_{G}(\{ v_{i}^{-}, v_{i}^{-2} \}; \{ v_{3-i}^{-}, v_{3-i}^{-2} \}) \neq \emptyset$, 
which contradicts Lemma \ref{maximality of |V(C)| (2)}(\ref{no v1-v2-}) 
or 
Lemma \ref{v1 v2 in P5-free}(\ref{no v1-2v2-2}). 
Thus $v_{1}v_{2}^{-} \notin E(G)$ for $i\in \{1,2\}$, and hence $v_{1}v_{2} \in E(G)$. 
\qed

\begin{Claim}
\label{v1v1-2}
For $i \in \{1, 2\}$, $v_{i}v_{i}^{-2} \in E(G)$. 
\end{Claim}
\proof 
Suppose that $v_{i}v_{i}^{-2} \notin E(G)$, 
and consider the paths $Q_{1} = v_{i}v_{i}^{-}v_{i}^{-2}$ and $Q_{2} = v_{i}v_{3-i}v_{3-i}^{-}$ 
(note that by Claim \ref{v1v2}, $v_{1}v_{2} \in E(G)$). 
Then $Q_{1}$ is an induced path. 
Hence by Lemma \ref{consecutive}, 
$N_{G}( \{v_{i}^{-}, v_{i}^{-2}\}; \{ v_{3-i}, v_{3-i}^{-} \}) \neq \emptyset$ or $v_{i}v_{3-i}^{-} \in E(G)$. 
This together with Lemma \ref{maximality of |V(C)| (2)}(\ref{no v1-v2-}) implies that 
$E(G) \cap \{v_{i}v_{3-i}^{-}, v_{i}^{-}v_{3-i}, v_{i}^{-2}v_{3-i} \} \neq \emptyset$. 

Assume first that 
$v_{i}^{-}v_{3-i} \in E(G)$ or $v_{i}^{-2}v_{3-i} \in E(G)$. 
Then by Claim \ref{no v1v2- or v1v2-2}, 
$G[\{v_{3-i}, v_{i}^{-}, v_{i}^{-2}\}]$ contains an induced path $Q_{1}'$ of order $3$ with an end $v_{3-i}$. 
By Claim \ref{no v1v1-2if v1v2- or v1v2-2}, 
we also see that $Q_{2}' = v_{3-i}v_{3-i}^{-}v_{3-i}^{-2}$ is an induced path. 
Hence by Lemma \ref{consecutive}, 
$N_{G}( \{v_{i}^{-}, v_{i}^{-2}\}; \{v_{3-i}^{-}, v_{3-i}^{-2}\}) \neq \emptyset$, 
which contradicts 
Lemma \ref{maximality of |V(C)| (2)}(\ref{no v1-v2-}) 
or 
Lemma \ref{v1 v2 in P5-free}(\ref{no v1-2v2-2}). 
Thus $E(G) \cap \{v_{i}^{-}v_{3-i}, v_{i}^{-2}v_{3-i}\} = \emptyset$, 
and hence $v_{i}v_{3-i}^{-} \in E(G)$. 
We now consider the paths $Q_{1}$ and $Q_{2}'' = v_{i}v_{3-i}^{-}v_{3-i}^{-2}$. 
Then by Claim \ref{no v1v2- or v1v2-2}, 
$Q_{2}''$ is an induced path. 
Hence by Lemma \ref{consecutive}, 
$N_{G}(\{ v_{i}^{-}, v_{i}^{-2} \}; \{v_{3-i}^{-}, v_{3-i}^{-2}\}) \neq \emptyset$, 
which contradicts Lemma \ref{maximality of |V(C)| (2)}(\ref{no v1-v2-}) 
or 
Lemma \ref{v1 v2 in P5-free}(\ref{no v1-2v2-2}) again. 
Thus $v_{i}v_{i}^{-2} \in E(G)$.  
\qed

Now we choose a longest cycle $\ora{C}$, a component $H$ 
and vertices $v_{1}$ and $v_{2}$ such that 
$|N_{G}(v_{1}; H) \cup N_{G}(v_{2}; H)| \ge 2$ and $N_{G}(H; v_{1}^{+} \ora{C} v_{2}^{-}) = \emptyset$ 
so that 
\begin{enumerate}
\renewcommand{\labelenumi}{\upshape{(C\arabic{enumi}})}
\setcounter{enumi}{\value{mymemory}}

\item
\label{choice 3} 
$|v_{1} \ora{C} v_{2}|$ is as large as possible, subject to (C\ref{choice 1}) and (C\ref{choice 2}).  
\end{enumerate}
Then by the choice, 
we can easily obtain the following.

\begin{Claim}
\label{no v1-v1+}
$v_{1}^{-}v_{1}^{+} \notin E(G)$. 
\end{Claim}
\proof 
Note that by Lemma \ref{maximality of |V(C)| (1)}, $N_{G}(H; \{v_{1}^{-}, v_{1}^{+}\}) = \emptyset$. 
If $v_{1}^{-}v_{1}^{+} \in E(G)$, 
then since $v_{1}v_{1}^{-2} \in E(G)$ by Claim \ref{v1v1-2}, 
$D = v_{1}^{-2}v_{1}v_{1}^{-}v_{1}^{+}\ora{C}v_{1}^{-2}$ is a cycle in $G$ 
such that $V(D) = V(C)$, 
$N_{G}(H; v_{1}^{+} \ora{D} v_{2}^{-}) = \emptyset$ and 
$|v_{1} \ora{D} v_{2}| > |v_{1} \ora{C} v_{2}|$, 
which contradicts the choice (C\ref{choice 3}).  
\qed

\begin{Claim}
\label{v1-2v1+}
$v_{1}^{-2}v_{1}^{+} \in E(G)$. 
\end{Claim}
\proof 
Let $xx' \in E(H)$. 
Note that by Claims \ref{xv1, xv2 in E(G)} and \ref{v1v1-2}, 
$\{xv_{1}, x'v_{1}, v_{1}v_{1}^{-2}\} \subseteq E(G)$. 
By Lemmas \ref{maximality of |V(C)| (1)} and \ref{N_{G}(H; C) cap N_{G}(H; C)^{-2} = emptyset}, 
$E(G) \cap \{xv_{1}^{-2}, xv_{1}^{-}, xv_{1}^{+}, x'v_{1}^{-2}, x'v_{1}^{-}, x'v_{1}^{+} \} = \emptyset$. 
By Claim \ref{no v1-v1+}, 
we also have $v_{1}^{-}v_{1}^{+} \notin E(G)$. 
Therefore, 
if $v_{1}^{-2}v_{1}^{+} \notin E(G)$, 
then 
$G[\{v_{1}, v_{1}^{-}, v_{1}^{-2}, x, x', v_{1}^{+}\}] \cong W^{*}$, a contradiction. 
\qed

Note that by Lemma \ref{maximality of |V(C)| (2)}(\ref{no v1-v2-}) 
and Claim \ref{v1-2v1+}, 
$|v_{2} \ora{C} v_{1}| \ge 6$. 
By Lemma \ref{maximality of |V(C)| (2)}(\ref{no v1-v2-}), we have 
\begin{align}
\label{v_{1}^{-2}v_{2}^{-} notin E(G)}
E(G) \cap \{v_{1}^{-}v_{2}^{-}, v_{1}^{-2}v_{2}^{-}\} = \emptyset. 
\end{align}
Since $v_{1}^{-2}v_{1}^{+} \in E(G)$ by Claim \ref{v1-2v1+}, 
it follows from Lemma \ref{maximality of |V(C)| (2)}(\ref{no v1-w v2+w+}) that 
\begin{align}
\label{v_{1}^{-3}v_{2}^{-} notin E(G)}
v_{1}^{-3}v_{2}^{-} \notin E(G). 
\end{align}
Since 
$v_{i}v_{i}^{-2} \in E(G)$ for $i \in \{1, 2\}$
by Claim \ref{v1v1-2}, 
it follows from Claim \ref{no v1v1-2if v1v2- or v1v2-2} that 
\begin{align}
\label{v_{1}v_{2}^{-}, v_{2}v_{1}^{-}, v_{2}v_{1}^{-2} notion E(G)}
E(G) \cap \{v_{1}v_{2}^{-}, v_{2}v_{1}^{-}, v_{2}v_{1}^{-2} \} = \emptyset. 
\end{align}

Consider the paths $Q_{1} = v_{1}v_{2}v_{2}^{-}$ and $Q_{2} = v_{1}v_{1}^{-2}v_{1}^{-3}$. 
By (\ref{v_{1}v_{2}^{-}, v_{2}v_{1}^{-}, v_{2}v_{1}^{-2} notion E(G)}), 
$Q_{1}$ is an induced path. 
Hence 
by Lemma \ref{consecutive}, 
$N_{G}(\{v_{2}, v_{2}^{-}\}; \{v_{1}^{-2}, v_{1}^{-3}\}) \neq \emptyset$ or $v_{1}v_{1}^{-3} \in E(G)$. 
This together with (\ref{v_{1}^{-2}v_{2}^{-} notin E(G)})--(\ref{v_{1}v_{2}^{-}, v_{2}v_{1}^{-}, v_{2}v_{1}^{-2} notion E(G)}) 
implies that 
$E(G) \cap \{v_{1}v_{1}^{-3}, v_{2}v_{1}^{-3}\} \neq \emptyset$. 
Note that 
by Lemmas \ref{maximality of |V(C)| (1)} and \ref{N_{G}(H; C) cap N_{G}(H; C)^{-2} = emptyset}, 
$N_{G}(H; \{v_{1}^{-2}, v_{1}^{-}\}) = \emptyset$. 
Note also that by Claim \ref{v1-2v1+}, $v_{1}^{-2}v_{1}^{+} \in E(G)$. 
Therefore, 
if $v_{1}^{-3}v_{1} \in E(G)$, 
then $D = v_{1}^{-3}v_{1} v_{1}^{-} v_{1}^{-2} v_{1}^{+} \ora{C} v_{1}^{-3}$ is a cycle in $G$ such that 
$V(D) = V(C)$, 
$N_{G}(H; v_{1}^{+} \ora{D} v_{2}^{-}) = \emptyset$ and 
$|v_{1} \ora{D} v_{2}| > |v_{1} \ora{C} v_{2}|$, 
which contradicts the choice (C\ref{choice 3}). 
Thus $v_{1}^{-3}v_{1} \notin E(G)$, and hence $v_{1}^{-3}v_{2} \in E(G)$.

We next consider 
the paths $Q_{1}' = v_{1}^{-3}v_{2}v_{2}^{-}$ and $Q_{2}' = v_{1}^{-3}v_{1}^{-2}v_{1}^{-}$. 
By (\ref{v_{1}^{-3}v_{2}^{-} notin E(G)}), 
$Q_{1}'$ is an induced path. 
Since $N_{G}(\{v_{2}, v_{2}^{-}\}; \{v_{1}^{-}, v_{1}^{-2} \}) = \emptyset$ 
by (\ref{v_{1}^{-2}v_{2}^{-} notin E(G)}) and (\ref{v_{1}v_{2}^{-}, v_{2}v_{1}^{-}, v_{2}v_{1}^{-2} notion E(G)}), 
it follows from Lemma \ref{consecutive} that $v_{1}^{-}v_{1}^{-3} \in E(G)$. 
Note that by Claims \ref{v1v1-2} and \ref{v1-2v1+}, 
$\{v_{1}v_{1}^{-2}, v_{1}^{-2}v_{1}^{+}\} \subseteq E(G)$, 
and 
hence $D = v_{1}^{-3}v_{1}^{-}v_{1}v_{1}^{-2}v_{1}^{+} \ora{C} v_{1}^{-3}$ is a cycle in $G$ 
such that $V(D) = V(C)$,  
$N_{G}(H; v_{1}^{+} \ora{D} v_{2}^{-}) = \emptyset$ and 
$|v_{1} \ora{D} v_{2}| > |v_{1} \ora{C} v_{2}|$, 
which contradicts the choice (C\ref{choice 3}). 

This completes the proof of Theorem \ref{P_{5}, W^{*}}. 
\qed

\section{Proof of Theorem \ref{P_{5}, K_{4}^{-}}}
\label{proof of P_{5},K_{4}^{-}}

Let $G$ be a $2$-connected $\{P_{5}, K_{4}^{-}\}$-free graph.
We first introduce a useful claim for our proof.

\begin{Claim}
\label{complete} 
Let $\ora{Q}$ be a path of $G$ starting from $v \in V(G)$ 
such that $v$ is adjacent to every vertex in $V(Q) \setminus \{v\}$, 
and let $a \in V(G) \setminus V(Q)$. 
Then either 
$V(Q) \subseteq N_{G}(a)$ 
or $|N_{G}(a; Q)| \le 1$. 
\end{Claim}
\proof
If $|Q| \le 2$, then the assertion clearly holds. Thus we may assume that $|Q| \ge 3$.  

We first suppose that $G[V(Q)]$ is not complete. 
Then there exist $h, l$ with $1 \le h < l \le |Q| - 1$ such that $v^{+h}v^{+l} \not\in E(G)$.
Choose $h$ and $l$ so that $l-h$ is as small as possible. 
Note that $l \geq h+2$ and $v^{+h}v^{+(h+1)}, v^{+(h+1)}v^{+l} \in E(G)$. 
Hence $\{v,v^{+h},v^{+(h+1)},v^{+l}\}$ induces $K_{4}^{-}$ in $G$, which is a contradiction.
Thus $G[V(Q)]$ is complete.

If $V(Q) \not\subseteq N_{G}(a)$ and $|N_{G}(a; Q)| \ge 2$, 
then there exist three vertices $u$, $u'$ and $u''$ such that $u,u'\in N_{G}(a; Q)$ and $u'' \notin N_{G}(a; Q)$,
and hence $\{a,u,u',u''\}$ induces $K_{4}^{-}$ in $G$ because $G[V(Q)]$ is complete, which is a contradiction.
Consequently, we get the desired conclusion.
\qed

We show that $G$ contains a longest cycle which is a dominating cycle. 
By way of a contradiction, 
suppose that any longest cycles of $G$ are not dominating cycles.
Let $\ora{C}$ be the same described as in 
the paragraph preceding Lemma \ref{indep. set S} in Subsection \ref{longest cycle P5-free no dom}, 
and let $H$ be a component of $G-C$ such that $|H| = \mu(C) \ ( \ \ge 2)$. 
Since $G$ is $2$-connected, 
there exist two distinct vertices $v_{1}$ and $v_{2}$ in $N_{G}(H; C)$ 
such that $|N_{G}(v_{1}; H) \cup N_{G}(v_{2}; H)| \ge 2$. 
Then 
$|v_{i} \ora{C} v_{3-i}| \ge 4$ for $i \in \{1, 2\}$ because $C$ is longest.  
(Note that by these assumptions, 
we can use all lemmas of Section \ref{preparation}.) 
We choose the vertices $v_{1}$ and $v_{2}$ so that 
\begin{align}
\label{N_{G}(H; v_{1}^{+} C v_{2}^{-}) is emptyset (2)}
N_{G}(H; v_{1}^{+} \ora{C} v_{2}^{-}) = \emptyset. 
\end{align}

\begin{Claim}
\label{|V(P)| = 2}
There exists an edge $x_{1}x_{2}$ in $H$ such that $v_{i}x_{i} \in E(G)$ for $i \in \{1, 2\}$. 
\end{Claim}
\proof 
Suppose not.
Let $x_{1}\in N_{G}(v_{1};H)$ and $x_{2}\in N_{G}(v_{2};H)$ be distinct vertices, 
and let $P$ be a shortest $(x_{1}, x_{2})$-path in $H$.
We choose $x_{1}$ and $x_{2}$ so that $|P|$ is as small as possible.
Then $x_{1}x_{2} \notin E(G)$ and 
$\emptyset \neq V(P) \setminus \{x_{1}, x_{2}\} \subseteq V(P) \setminus N_{G}(v_{i})$ for $i \in \{1, 2\}$. 
Hence 
by Lemma \ref{if V(H) - N_{G}(v_{i}) neq emptyset} and (\ref{N_{G}(H; v_{1}^{+} C v_{2}^{-}) is emptyset (2)}), 
$V(v_{1}^{+} \ora{C} v_{2}^{-}) \subseteq N_{G}(v_{1}) \cap N_{G}(v_{2})$, 
and this implies that $v_{1} v_{2} \in E(G)$ 
(otherwise, $G[V(v_{1} \ora{C} v_{2})]$ contains a $K_{4}^{-}$ as an induced subgraph, a contradiction). 
On the other hand, 
consider the paths $Q_{1} = P$ and $Q_{2} = x_{1}v_{1}v_{1}^{+}$. 
Then by the minimality of $|P|$, 
$Q_{1}$ is an induced path of order at least $3$. 
By Lemma \ref{maximality of |V(C)| (1)}, 
$Q_{2}$ is also an induced path. 
Hence by Lemma \ref{consecutive}, 
$N_{G}( P - x_{1}; \{v_{1}, v_{1}^{+}\}) \neq \emptyset$. 
Combining this with (\ref{N_{G}(H; v_{1}^{+} C v_{2}^{-}) is emptyset (2)}) 
and the fact that $V(P) \setminus \{x_{1}, x_{2}\} \subseteq V(P) \setminus N_{G}(v_{1})$, 
we get $v_{1}x_{2} \in E(G)$. 
Similarly, 
by considering the paths $P$ and $x_{2}v_{2}v_{2}^{-}$, 
we also have $v_{2}x_{1} \in E(G)$. 
But then $G[\{v_{1}, v_{2}, x_{1}, x_{2}\}] \cong K_{4}^{-}$, 
a contradiction.  
\qed

Let $x_{1}x_{2}$ be as in Claim \ref{|V(P)| = 2}. 
By the symmetry of $\ora{C}$ and $\ola{C}$, 
we may always assume that $v_{1}x_{2} \notin E(G)$ 
if $\{v_{1} x_{2}, v_{2} x_{1}\} \not \subseteq E(G)$.

Now let $w_{1}$ be a vertex in $v_{2}^{+} \ora{C} v_{1}^{-}$ 
such that 
$V(w_{1}\ora{C}v_{1}^{-})\subseteq N_{G}(v_{1})$.
We choose $w_{1}$ so that $|w_{1} \ora{C} v_{1}|$ is as large as possible. 
By Lemma \ref{maximality of |V(C)| (2)}(\ref{no v1-v2-}), Claim \ref{complete} 
and the choice of $w_{1}$, we can easily obtain the following.

\begin{figure}[htbp]
\begin{center}
\hspace{-10pt}\includegraphics[scale=0.9,clip]{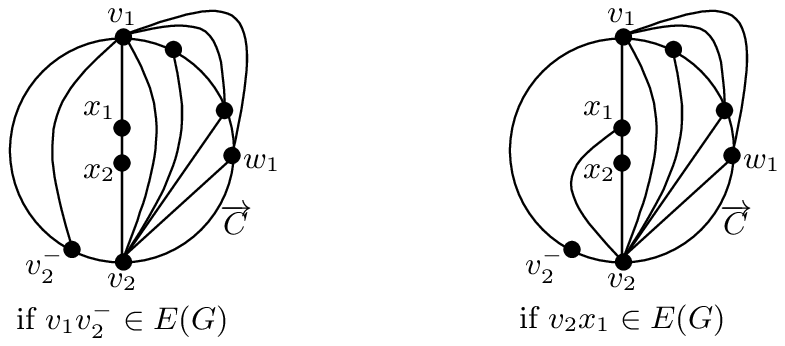}
\vspace{-25pt}
\caption{Claim \ref{easy fact}(\ref{|N_{G}(v_{2}; w_{1} C v_{1})| <= 1})}
\label{w1(Fig)}
\end{center}
\end{figure}

\begin{Claim}
\label{easy fact}
\begin{enumerate}[{\upshape(i)}]

\item
\label{N_{G}(x1; w1 C v1) = v1} 
$|N_{G}(x; w_{1} \ora{C}v_{1})| \le 1$ for $x \in \{v_{2}^{-}, v_{2}^{-2}, x_{1}, x_{2}\}$.

\item
\label{N_{G}(w1-; w1 C v1) = w1} 
If $w_{1}^{-} \neq v_{2}$, then $N_{G}(w_{1}^{-}; w_{1} \ora{C} v_{1}) = \{w_{1}\}$.

\item 
\label{|N_{G}(v_{2}; w_{1} C v_{1})| <= 1}
If 
$\{v_{1}v_{2}^{-}, v_{2}x_{1}\} \cap E(G) \neq \emptyset$, 
then $|N_{G}(v_{2}; w_{1} \ora{C} v_{1})| \le 1$. 

\end{enumerate}
\end{Claim}
\proof 
Let $x \in \{v_{2}^{-}, v_{2}^{-2}, x_{1}, x_{2}\}$. 
Then 
by Lemmas \ref{maximality of |V(C)| (1)} and 
\ref{maximality of |V(C)| (2)}(\ref{no v1-v2-}), 
$xv_{1}^{-} \notin E(G)$. 
Hence by applying Claim \ref{complete} as $\ora{Q} = v_{1}\ola{C}w_{1}$ and $a = x$, 
we have $|N_{G}(x; w_{1} \ora{C}v_{1})| \le 1$. 
Thus (\ref{N_{G}(x1; w1 C v1) = v1}) holds. 
If $w_{1}^{-} \neq v_{2}$, 
then by the choice of $w_{1}$, 
$w_{1}^{-}v_{1} \notin E(G)$, 
and hence again by Claim \ref{complete}, 
$N_{G}(w_{1}^{-}; w_{1} \ora{C} v_{1}) = \{w_{1}\}$. 
Thus (\ref{N_{G}(w1-; w1 C v1) = w1}) also holds. 
To show (\ref{|N_{G}(v_{2}; w_{1} C v_{1})| <= 1}), 
suppose that 
$\{v_{1}v_{2}^{-}, v_{2}x_{1}\} \cap E(G) \neq \emptyset$ 
and 
$|N_{G}(v_{2}; w_{1} \ora{C} v_{1})| \ge 2$. 
Since $|N_{G}(v_{2}; w_{1} \ora{C} v_{1})| \ge 2$, 
it follows from Claim \ref{complete} that 
$V(w_{1} \ora{C} v_{1}) \subseteq N_{G}(v_{2})$, 
and thus $v_{2}w_{1} \ora{C} v_{1}$ is a path with an end $v_{2}$ such that $V(w_{1} \ora{C} v_{1}) \subseteq N_{G}(v_{2})$ 
(see Figure \ref{w1(Fig)}). 
Since 
$\{v_{1}x_{1}, v_{2}v_{2}^{-} \} \subseteq E(G)$, 
the assumption $\{v_{1}v_{2}^{-}, v_{2}x_{1}\} \cap E(G) \neq \emptyset$ 
implies that 
$\{v_{1}, v_{2}\} \subseteq N_{G}(x)$ for some $x \in \{x_{1}, v_{2}^{-}\}$, 
and thus $|N_{G}(x; v_{2}w_{1} \ora{C} v_{1} )| \ge |\{v_{1}, v_{2}\}| = 2$. 
Then again by Claim \ref{complete}, 
$V(v_{2}w_{1} \ora{C} v_{1}) \subseteq N_{G}(x)$, 
in particular, $xv_{1}^{-} \in E(G)$, which contradicts 
Lemma \ref{maximality of |V(C)| (1)} or Lemma \ref{maximality of |V(C)| (2)}(\ref{no v1-v2-}). 
Thus (\ref{|N_{G}(v_{2}; w_{1} C v_{1})| <= 1}) holds.  
\qed

Since $v_{1}x_{1} \in E(G)$, 
the following fact is directly obtained from Claim \ref{easy fact}(\ref{N_{G}(x1; w1 C v1) = v1}).

\begin{Fact}
\label{N_{G}( x_{1}, v_{2}^{-} ; w_{1}, w_{1}^{+}) = emptyset}
$N_{G}(x_{1}; w_{1} \ora{C} v_{1}^{-}) = \emptyset$.  
\end{Fact}

We divide the proof into 
two cases according as $\{v_{1}x_{2}, v_{2}x_{1}\} \not \subseteq E(G)$ or $\{v_{1}x_{2}, v_{2}x_{1}\} \subseteq E(G)$.

\bigskip
\noindent
\textbf{Case 1.} $\{v_{1}x_{2}, v_{2}x_{1}\} \not \subseteq E(G)$. 

Then $v_{1} x_{2} \notin E(G)$ (see the paragraph following the proof of Claim \ref{|V(P)| = 2}). 
Hence by Lemma \ref{if V(H) - N_{G}(v_{i}) neq emptyset} and (\ref{N_{G}(H; v_{1}^{+} C v_{2}^{-}) is emptyset (2)}),
\begin{align}
\label{v2-v1 v2-2v1}
\textup{$\{v_{2}^{-}, v_{2}^{-2}\} \subseteq N_{G}(v_{1})$.}
\end{align} 
Then by applying Claim \ref{easy fact}(\ref{N_{G}(x1; w1 C v1) = v1}) as $x = v_{2}^{-}$ and $x = v_{2}^{-2}$, 
the following fact holds.

\begin{Fact}
\label{N_{G}( v_{2}^{-2}, v_{2}^{-} ; w_{1}, w_{1}^{+}) = emptyset}
$N_{G}(v_{2}^{-h}; w_{1} \ora{C} v_{1}^{-}) = \emptyset$ for $h \in \{1, 2\}$.  
\end{Fact}

Moreover, 
by Lemmas \ref{maximality of |V(C)| (1)} and \ref{N_{G}(H; C) cap N_{G}(H; C)^{-2} = emptyset}, 
$N_{G}(H; \{v_{1}^{-}, v_{1}^{-2}\}) = \emptyset$, 
and hence 
Lemma \ref{if V(H) - N_{G}(v_{i}) neq emptyset} yields that $v_{1}v_{1}^{-2} \in E(G)$. 
This together with the choice of $w_{1}$ implies that 
\begin{align}
\label{w1+ neq v1}
w_{1}^{+} \neq v_{1}. 
\end{align}

\begin{figure}[htbp]
\begin{center}
\hspace{-10pt}\includegraphics[scale=0.9,clip]{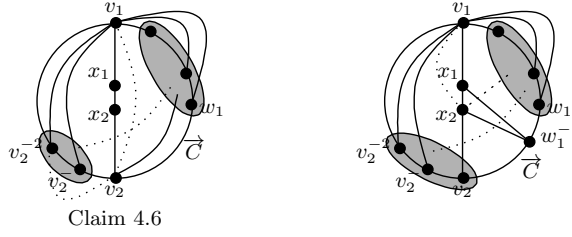}
\vspace{-25pt}
\caption{The cycle $C$ in Case 1}
\label{case1(Fig)}
\end{center}
\end{figure}

\begin{Claim}
\label{N_{G}(v_{2}; w_{1} C v_{1}^{-}) = emptyset}
$N_{G}(v_{2}; w_{1}\ora{C} v_{1}^{-}) = \emptyset$. In particular, $w_{1}^{-} \neq v_{2}$. 
\end{Claim}
\proof 
Suppose that 
$N_{G}(v_{2}; w_{1}\ora{C} v_{1}^{-}) \neq \emptyset$. 
Then by Claim \ref{easy fact}(\ref{|N_{G}(v_{2}; w_{1} C v_{1})| <= 1}) and (\ref{v2-v1 v2-2v1}), 
$|N_{G}(v_{2}; w_{1} \ora{C} v_{1})| = |N_{G}(v_{2}; w_{1} \ora{C} v_{1}^{-})| = 1$. 
By (\ref{w1+ neq v1}), 
the equality $|N_{G}(v_{2}; w_{1} \ora{C} v_{1}^{-})| = 1$ implies that 
$G[V(w_{1} \ora{C} v_{1}^{-}) \cup \{v_{2}\}]$ contains an induced path $Q_{1}$ of order at least $3$ 
with an end $v_{2}$. 
On the other hand, 
the equality $|N_{G}(v_{2}; w_{1} \ora{C} v_{1})| = |N_{G}(v_{2}; w_{1} \ora{C} v_{1}^{-})|$ implies that 
$v_{1}v_{2} \notin E(G)$. 
Since 
$v_{2}v_{2}^{-} \in E(G)$ 
and 
$G[\{v_{1}, v_{2}^{-}, v_{2}^{-2}\}]$ is triangle by (\ref{v2-v1 v2-2v1}), 
these together with Claim \ref{complete} imply that $v_{2}v_{2}^{-2} \notin E(G)$ 
(see the left of Figure \ref{case1(Fig)}), 
and thus 
$Q_{2} = v_{2}v_{2}^{-}v_{2}^{-2}$ is also an induced path. 
Hence 
by Lemma \ref{consecutive}, 
$N_{G}(Q_{1} - v_{2}; Q_{2} - v_{2}) \neq \emptyset$. 
Since $N_{G}(Q_{1} - v_{2}; Q_{2} - v_{2}) \subseteq N_{G}(\{v_{2}^{-}, v_{2}^{-2} \}; w_{1}\ora{C}v_{1}^{-})$, 
this contradicts 
Fact \ref{N_{G}( v_{2}^{-2}, v_{2}^{-} ; w_{1}, w_{1}^{+}) = emptyset}. 
\qed

\begin{Claim}
\label{N_{G}(x_{2}; w_{1} C v_{1}^{-}) = emptyset}
$N_{G}(x_{2}; w_{1}\ora{C} v_{1}^{-}) = \emptyset$. 
\end{Claim}
\proof 
Suppose not. 
Then by Claim \ref{easy fact}(\ref{N_{G}(x1; w1 C v1) = v1}), 
$|N_{G}(x_{2}; w_{1}\ora{C} v_{1}^{-})| = 1$, 
and this together with 
(\ref{w1+ neq v1}) implies that 
$G[V(w_{1} \ora{C} v_{1}^{-}) \cup \{x_{2}\}]$ contains an induced path $Q_{1}$ of order at least $3$ 
with an end $x_{2}$. 
On the other hand, 
by (\ref{N_{G}(H; v_{1}^{+} C v_{2}^{-}) is emptyset (2)}), 
the path $Q_{2} = x_{2}v_{2}v_{2}^{-}$ is also an induced path. 
Therefore, by Lemma \ref{consecutive}, 
we have $N_{G}(\{v_{2}, v_{2}^{-}\}; w_{1}\ora{C} v_{1}^{-}) \supseteq N_{G}(Q_{1}-x_{2}; Q_{2} - x_{2}) \neq \emptyset$, 
which contradicts Fact \ref{N_{G}( v_{2}^{-2}, v_{2}^{-} ; w_{1}, w_{1}^{+}) = emptyset} 
or Claim \ref{N_{G}(v_{2}; w_{1} C v_{1}^{-}) = emptyset}. 
\qed

\begin{Claim}
\label{x_{1}w_{1}^{-}, x_{2}w_{2}^{-}}
$\{x_{1}w_{1}^{-}, x_{2}w_{1}^{-}\} \subseteq E(G)$. 
\end{Claim}
\proof 
Consider the paths $Q_{1} = v_{1}w_{1}w_{1}^{-}$ and $Q_{2} = v_{1}x_{1}x_{2}$. 
Since $w_{1}^{-} \neq v_{2}$ by Claim \ref{N_{G}(v_{2}; w_{1} C v_{1}^{-}) = emptyset}, 
it follows from Claim \ref{easy fact}(\ref{N_{G}(w1-; w1 C v1) = w1}) that 
$Q_{1}$ is an induced path. 
Since $v_{1}x_{2} \notin E(G)$, 
$Q_{2}$ is also an induced path. 
Hence by Lemma \ref{consecutive}, 
Fact \ref{N_{G}( x_{1}, v_{2}^{-} ; w_{1}, w_{1}^{+}) = emptyset} and Claim \ref{N_{G}(x_{2}; w_{1} C v_{1}^{-}) = emptyset}, 
it follows that $N_{G}(w_{1}^{-}; \{x_{1}, x_{2}\}) \neq \emptyset$. 
We next consider the path $w_{1}^{+}w_{1}w_{1}^{-}$ and a path in $G[\{w_{1}^{-}, x_{1}, x_{2}\}]$ of order $3$ with an end $w_{1}^{-}$.
By Claim \ref{easy fact}(\ref{N_{G}(w1-; w1 C v1) = w1}), the path $w_{1}^{+}w_{1}w_{1}^{-}$ is an induced path.
Hence by again Lemma \ref{consecutive}, Fact \ref{N_{G}( x_{1}, v_{2}^{-} ; w_{1}, w_{1}^{+}) = emptyset} and Claim \ref{N_{G}(x_{2}; w_{1} C v_{1}^{-}) = emptyset},
we can easily obtain 
$\{x_{1}w_{1}^{-}, x_{2}w_{2}^{-}\} \subseteq E(G)$. 
\qed

The graph illustrated in the right of Figure \ref{case1(Fig)} is 
a current situation. 
We now consider the paths $Q_{1} = w_{1}^{-}w_{1}w_{1}^{+}$ 
and $Q_{2} = w_{1}^{-}x_{2}v_{2}$. 
By Claims \ref{easy fact}(\ref{N_{G}(w1-; w1 C v1) = w1}) and \ref{N_{G}(v_{2}; w_{1} C v_{1}^{-}) = emptyset}, 
$Q_{1}$ is an induced path. 
Hence by Lemma \ref{consecutive} 
and Claims \ref{N_{G}(v_{2}; w_{1} C v_{1}^{-}) = emptyset} and \ref{N_{G}(x_{2}; w_{1} C v_{1}^{-}) = emptyset}, 
we have $w_{1}^{-}v_{2} \in E(G)$, 
and hence $|N_{G}(v_{2}; \{x_{1}, x_{2}, w_{1}^{-}\})| \ge |\{x_{2}, w_{1}^{-}\}| = 2$. 
Since $G[\{x_{1}, x_{2}, w_{1}^{-}\}]$ is a triangle by Claim \ref{x_{1}w_{1}^{-}, x_{2}w_{2}^{-}}, 
it follows from Claim \ref{complete} that 
$v_{2}x_{1} \in E(G)$. 
Therefore, we see that $G[\{v_{2}, x_{1}, x_{2}\}]$ is also a triangle. 
Since 
$\{v_{2}, x_{1}, x_{2}\} \not \subseteq N_{G}(v_{1})$
because $v_{1}x_{2} \notin E(G)$, 
Claim \ref{complete} also yields that 
$v_{1}v_{2}$ is not an edge in $G$.

On the other hand, 
consider the paths $Q_{1}' = v_{2}^{-}v_{2}x_{2}$ and $Q_{2}' = v_{2}^{-}v_{1}w_{1}$. 
Since $v_{2}^{-}x_{2} \notin E(G)$ by (\ref{N_{G}(H; v_{1}^{+} C v_{2}^{-}) is emptyset (2)}), 
it follows that $Q_{1}'$ is an induced path. 
By Fact \ref{N_{G}( v_{2}^{-2}, v_{2}^{-} ; w_{1}, w_{1}^{+}) = emptyset}, 
$Q_{2}'$ is also an induced path. 
Hence by Lemma \ref{consecutive}, 
$N_{G}(\{v_{2}, x_{2}\}; \{v_{1}, w_{1}\}) \neq \emptyset$. 
This together with 
Claims \ref{N_{G}(v_{2}; w_{1} C v_{1}^{-}) = emptyset}, \ref{N_{G}(x_{2}; w_{1} C v_{1}^{-}) = emptyset} 
and the assumption $v_{1}x_{2} \notin E(G)$ implies 
that 
$v_{1}v_{2}$ is an edge in $G$. 
This is a contradiction.

\bigskip
\noindent
\textbf{Case 2.} $\{v_{1}x_{2}, v_{2}x_{1}\} \subseteq E(G)$. 

By the assumption of Case 2, 
Claim \ref{easy fact}(\ref{N_{G}(w1-; w1 C v1) = w1}) and (\ref{|N_{G}(v_{2}; w_{1} C v_{1})| <= 1}), 
the following claim holds. 

\begin{Claim}
\label{easy fact case 2}
$N_{G}(v_{2}; w_{1} \ora{C} v_{1}) = \{v_{1}\}$ and $N_{G}(w_{1}^{-}; w_{1} \ora{C} v_{1}) = \{w_{1}\}$.
\end{Claim}
\proof 
If $v_{1}v_{2}\not\in E(G)$, then $G[\{v_{1}, v_{2}, x_{1}, x_{2}\}]$ 
is isomorphic to $K_{4}^{-}$, a contradiction.
Thus $v_{1}v_{2}\in E(G)$.
This together with Claim \ref{easy fact}(\ref{|N_{G}(v_{2}; w_{1} C v_{1})| <= 1}) 
forces $N_{G}(v_{2}; w_{1} \ora{C} v_{1}) = \{v_{1}\}$.
In particular, $w_{1}^{-} \neq v_{2}$.
Hence by Claim \ref{easy fact}(\ref{N_{G}(w1-; w1 C v1) = w1}), 
we have $N_{G}(w_{1}^{-}; w_{1} \ora{C} v_{1}) = \{w_{1}\}$.
\qed

By Claim \ref{easy fact case 2}, 
$w_{1}^{-}v_{1} \notin E(G)$. 
Hence by applying Claim \ref{complete} as $Q = v_{1}x_{1}x_{2}$ and $a = w_{1}^{-}$, 
it follows that 
$|N_{G}(w_{1}^{-}; \{ v_{1}, x_{1}, x_{2} \})| \le 1$. 
Therefore, by changing the label of $x_{1}$ and $x_{2}$ if necessary, 
we may assume that 
\begin{align}
\label{no w1-x1}
w_{1}^{-}x_{1} \notin E(G). 
\end{align}

Let $z_{1}$ be a vertex in $v_{1}^{+} \ora{C} v_{2}^{-}$ 
such that 
$v_{1}z_{1} \in E(G)$. 
Then 
\begin{align}
\label{z1 neq v2-}
z_{1} \neq v_{2}^{-}
\end{align}
(otherwise, by (\ref{N_{G}(H; v_{1}^{+} C v_{2}^{-}) is emptyset (2)}), Claim \ref{easy fact case 2} and the assumption of Case 2, 
$G[\{v_{1}, v_{2}, x_{1}, z_{1}\}]$ is isomorphic to $K_{4}^{-}$, a contradiction). 
We choose $z_{1}$ so that $|v_{1} \ora{C} z_{1}|$ is as large as possible.

\begin{Claim}
\label{z2}
\begin{enumerate}[{\upshape(i)}]
\item
\label{V(z_{1}^{+} C v_{2}^{-}) subseteq N_{G}(z_{1})} 
$V(z_{1}^{+} \ora{C} v_{2}^{-}) \subseteq N_{G}(z_{1})$.

\item
\label{N_{G}(v_{2}; z_{1} C v_{2}^{-}) = v_{2}^{-}}
$N_{G}(v_{2}; z_{1} \ora{C}v_{2}^{-}) = \{v_{2}^{-}\}$. 

\end{enumerate}
\end{Claim}
\proof 
(\ref{V(z_{1}^{+} C v_{2}^{-}) subseteq N_{G}(z_{1})}) 
If $z_{1}^{+} = v_{2}^{-}$, 
then the assertion clearly holds. 
Thus we may assume that $z_{1}^{+} \neq v_{2}^{-}$. 
Consider the paths $Q_{1} = z_{1}v_{1}x_{1}$ and $Q_{2} = z_{1} \ora{C} v_{2}^{-}$. 
Then by (\ref{N_{G}(H; v_{1}^{+} C v_{2}^{-}) is emptyset (2)}), 
$Q_{1}$ is an induced path. 
By the choice of $z_{1}$ and again (\ref{N_{G}(H; v_{1}^{+} C v_{2}^{-}) is emptyset (2)}), 
$N_{G}(Q_{1} - z_{1}; Q_{2} - z_{1}) = \emptyset$. 
Hence Lemma \ref{consecutive} yields that 
$V(z_{1}^{+} \ora{C} v_{2}^{-}) = V(Q_{2} - z_{1}) \subseteq N_{G}(z_{1})$. 

(\ref{N_{G}(v_{2}; z_{1} C v_{2}^{-}) = v_{2}^{-}}) 
By Claim \ref{easy fact case 2}, (\ref{z1 neq v2-}) and the choice of $z_{1}$, 
we see that 
$|z_{1} \ora{C} v_{2}| - 1 \ge |N_{G}(v_{1}; v_{2}\ola{C}z_{1})| \ge |\{v_{2}, z_{1}\}| = 2$. 
Since $v_{2}\ola{C}z_{1}$ is a path with an end $v_{2}$, $V(v_{2}\ola{C}z_{1})\not\subseteq N_{G}(v_{1})$ and $|N_{G}(v_{1};v_{2}\ola{C}z_{1})|\geq 2$, 
it follows from Claim \ref{complete} that $V(z_{1}\ora{C}v_{2}^{-}) \not\subseteq N_{G}(v_{2})$, 
i.e., 
$|N_{G}(v_{2}; z_{1} \ora{C} v_{2}^{-})| \le |z_{1} \ora{C} v_{2}^{-}| - 1$. 
By Claim \ref{z2}(\ref{V(z_{1}^{+} C v_{2}^{-}) subseteq N_{G}(z_{1})}), 
$z_{1} \ora{C} v_{2}^{-}$ is a path with an end $z_{1}$ such that 
$V(z_{1}^{+} \ora{C} v_{2}^{-}) \subseteq N_{G}(z_{1})$, 
and hence Claim \ref{complete} implies that 
$N_{G}(v_{2}; z_{1} \ora{C}v_{2}^{-}) = \{v_{2}^{-}\}$. 
\qed

\begin{Claim}
\label{no v2-w1}
$N_{G}(v_{2}^{-h}; w_{1}^{-} \ora{C} v_{1}^{-}) = \emptyset$ for $h \in \{1, 2\}$. 
\end{Claim}
\proof
Suppose that $N_{G}(v_{2}^{-h}; w_{1}^{-} \ora{C} v_{1}^{-}) \neq \emptyset$ for some $h \in \{1, 2\}$, 
and choose $v_{2}^{-h}$ so that $h = 1$ if possible. 
Note that by Lemmas \ref{maximality of |V(C)| (2)}(\ref{no v1-v2-}) and \ref{v1 v2 in P5-free}(\ref{no v1-2v2-2}), 
$w_{1} \neq v_{1}^{-}$. 

\begin{figure}[htbp]
\begin{center}
\hspace{-108pt}\includegraphics[scale=0.9,clip]{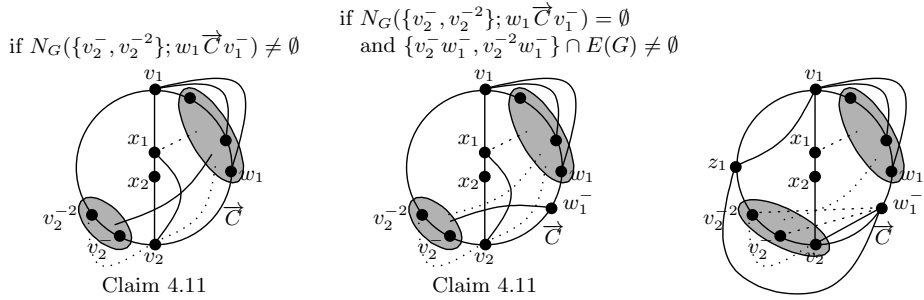}
\vspace{-25pt}
\caption{The cycle $C$ in Case 2}
\label{case2(Fig)}
\end{center}
\end{figure}

We first assume that $N_{G}(v_{2}^{-h}; w_{1} \ora{C} v_{1}^{-}) \neq \emptyset$.
Then by Claim \ref{easy fact}(\ref{N_{G}(x1; w1 C v1) = v1}), 
$|N_{G}(v_{2}^{-h}; w_{1} \ora{C} v_{1}^{-})| = 1$, 
and this implies that 
$G[V(w_{1} \ora{C} v_{1}^{-}) \cup \{v_{2}^{-h}\}]$ contains an induced path $Q_{1}$ of order at least $3$ 
with an end $v_{2}^{-h}$. 
By Claim \ref{z2}(\ref{N_{G}(v_{2}; z_{1} C v_{2}^{-}) = v_{2}^{-}}), 
(\ref{N_{G}(H; v_{1}^{+} C v_{2}^{-}) is emptyset (2)}) and (\ref{z1 neq v2-}), 
we also see that 
$Q_{2} = v_{2}^{-h} \ora{C}v_{2}x_{1}$ is an induced path of order at least $3$ (see the left of Figure \ref{case2(Fig)}). 
Hence by Lemma \ref{consecutive}, 
$N_{G}(Q_{1} - v_{2}^{-h}; Q_{2} - v_{2}^{-h}) \neq \emptyset$, 
which contradicts Fact \ref{N_{G}( x_{1}, v_{2}^{-} ; w_{1}, w_{1}^{+}) = emptyset}, 
Claim \ref{easy fact case 2} 
or the choice of $v_{2}^{-h}$. 
Thus $N_{G}(v_{2}^{-h}; w_{1} \ora{C} v_{1}^{-}) = \emptyset$, and so $v_{2}^{-h}w_{1}^{-}\in E(G)$.

Consider the paths $Q'_{1} = w_{1}^{-}w_{1}w_{1}^{+}$ 
and $Q'_{2} = w_{1}^{-}v_{2}^{-h}v_{2}^{-(h-1)}$. 
By Claim \ref{easy fact case 2}, 
$Q'_{1}$ is an induced path. 
Recall that $N_{G}(v_{2}^{-i}; w_{1} \ora{C} v_{1}^{-}) = \emptyset$ for $i\in \{1,2\}$.
This together with Claim \ref{easy fact case 2}
 leads to $N_{G}(Q'_{1} - w_{1}^{-}; Q'_{2} - w_{1}^{-}) = \emptyset$ 
(see also the center of Figure \ref{case2(Fig)}). 
Hence by Lemma \ref{consecutive}, 
we have $v_{2}^{-(h-1)}w_{1}^{-} \in E(G)$, 
and the choice of $v_{2}^{-h}$ implies that $h = 1$ and $v_{2}w_{1}^{-}\in E(G)$. 
Then, 
again by applying Lemma \ref{consecutive} as $(a, Q_{1}, Q_{2}) = (w_{1}^{-}, Q'_{1}, w_{1}^{-}v_{2}x_{1})$, 
we have $N_{G}(\{w_{1}, w_{1}^{+}\}; \{v_{2}, x_{1}\}) \neq \emptyset$ or $x_{1}w_{1}^{-} \in E(G)$. 
Then Fact \ref{N_{G}( x_{1}, v_{2}^{-} ; w_{1}, w_{1}^{+}) = emptyset} and Claim \ref{easy fact case 2} 
yield that $x_{1}w_{1}^{-} \in E(G)$,
which contradicts (\ref{no w1-x1}).
\qed

\begin{Claim}
\label{no z1w1}
$z_{1}w_{1} \notin E(G)$. 
\end{Claim}
\proof 
Suppose that $z_{1}w_{1} \in E(G)$, 
and consider the paths $Q_{1} = v_{2}^{-}z_{1}w_{1}$ 
(note that by Claim \ref{z2}(\ref{V(z_{1}^{+} C v_{2}^{-}) subseteq N_{G}(z_{1})}) and (\ref{z1 neq v2-}), $z_{1}v_{2}^{-} \in E(G)$) 
and $Q_{2} = v_{2}^{-}v_{2}x_{1}$. 
By Claim \ref{no v2-w1}, $Q_{1}$ is an induced path in $G$. 
By (\ref{N_{G}(H; v_{1}^{+} C v_{2}^{-}) is emptyset (2)}), 
$Q_{2}$ is also an induced path. 
Hence by Lemma \ref{consecutive}, 
$N_{G}(Q_{1} - v_{2}^{-}; Q_{2} - v_{2}^{-}) \neq \emptyset$, 
which contradicts 
Fact \ref{N_{G}( x_{1}, v_{2}^{-} ; w_{1}, w_{1}^{+}) = emptyset}, 
Claim \ref{easy fact case 2}, 
Claim \ref{z2}(\ref{N_{G}(v_{2}; z_{1} C v_{2}^{-}) = v_{2}^{-}}) 
or (\ref{N_{G}(H; v_{1}^{+} C v_{2}^{-}) is emptyset (2)}). 
\qed

\begin{Claim}
\label{z1w1-}
$\{z_{1}w_{1}^{-}, v_{2}w_{1}^{-}\} \subseteq E(G)$. 
\end{Claim}
\proof
We first show that $z_{1}w_{1}^{-} \in E(G)$. 
We consider the paths 
$Q_{1} = v_{1}w_{1}w_{1}^{-}$ 
and 
$Q_{2} = v_{1}z_{1}v_{2}^{-}$ 
(note that by Claim \ref{z2}(\ref{N_{G}(v_{2}; z_{1} C v_{2}^{-}) = v_{2}^{-}}) and (\ref{z1 neq v2-}), 
$z_{1}v_{2}^{-} \in E(G)$).  
By Claim \ref{easy fact case 2}, 
$Q_{1}$ is an induced path. 
By (\ref{z1 neq v2-}) and the choice of $z_{1}$, 
$Q_{2}$ is also an induced path. 
Hence by Lemma \ref{consecutive}, 
$N_{G}(\{z_{1}, v_{2}^{-}\}; \{w_{1}, w_{1}^{-}\}) \neq \emptyset$. 
Combining this with Claims \ref{no v2-w1} and \ref{no z1w1}, 
we get $z_{1}w_{1}^{-} \in E(G)$.

To show that $v_{2}w_{1}^{-} \in E(G)$, 
consider the paths $Q_{1}' = z_{1}v_{2}^{-}v_{2}$ and $Q_{2}' = z_{1}w_{1}^{-}w_{1}$. 
By Claim \ref{z2}(\ref{N_{G}(v_{2}; z_{1} C v_{2}^{-}) = v_{2}^{-}}) and (\ref{z1 neq v2-}), 
$Q_{1}'$ is an induced path of order $3$. 
By Claim \ref{no z1w1}, 
$Q_{2}'$ is also an induced path. 
Hence by Lemma \ref{consecutive}, 
$N_{G}(\{v_{2}, v_{2}^{-}\}; \{w_{1}, w_{1}^{-}\}) \neq \emptyset$. 
This together with
Claims \ref{easy fact case 2} and \ref{no v2-w1}
implies that 
$v_{2}w_{1}^{-} \in E(G)$. 
\qed

The graph illustrated in the right of Figure \ref{case2(Fig)} is 
a current situation. 
By Claim \ref{z1w1-}, $z_{1}w_{1}^{-}\in E(G)$.
This together with Claim \ref{no v2-w1} implies that $z_{1}\not\in \{v_{2}^{-},v_{2}^{-2}\}$.
In particular, $|z_{1} \ora{C} v_{2}| \ge 4$.
Now we consider the paths $Q_{1} = v_{2}w_{1}^{-}w_{1}$ 
(note that by Claim \ref{z1w1-}, $v_{2}w_{1}^{-} \in E(G)$) 
and $Q_{2} = v_{2}v_{2}^{-}v_{2}^{-2}$. 
By Claim \ref{easy fact case 2}, 
$Q_{1}$ is an induced path. 
By Claim \ref{z2}(\ref{N_{G}(v_{2}; z_{1} C v_{2}^{-}) = v_{2}^{-}}), 
$Q_{2}$ is also an induced path. 
Hence by Lemma \ref{consecutive}, 
$N_{G}(\{w_{1}^{-}, w_{1}\}; \{v_{2}^{-}, v_{2}^{-2}\})\not=\emptyset $, 
which contradicts Claim \ref{no v2-w1}.

This completes the proof of Theorem \ref{P_{5}, K_{4}^{-}}. 
\qed



\end{document}